\documentclass[11pt]{article}
\usepackage{amsfonts}
\usepackage[mathscr]{eucal}
\usepackage{mathrsfs}
\usepackage{amsmath}
\usepackage{amssymb}

\newtheorem{theorem}{Theorem}[section]
\newtheorem{lemma}[theorem]{Lemma}
\newtheorem{proposition}[theorem]{Proposition}
\newtheorem{definition}[theorem]{Definition}
\newtheorem{remark}[theorem]{Remark}
\newtheorem{corollary}[theorem]{Corollary}

\newtheorem{assumption}[theorem]{Assumption}
\newcommand{\De}{\Delta}

\newcommand{\Om}{{\Omega}}

\newcommand{\al}{{\alpha}}
\newcommand{\e}{\epsilon}
\newcommand{\ga}{{\gamma}}
\newcommand{\eps}{\varepsilon}
\newcommand{\la}{\lambda}
\newcommand{\si}{{\sigma}}
\newcommand{\tht}{{\theta}}

\newcommand{\Ac}{\mathfrak{A}}

\newcommand{\sC}{\mathscr{C}}
\newcommand{\sD}{\mathscr{D}}

\newcommand{\sQ}{\mathscr{Q}}

\newcommand{\N}{{\mathbb N}}

\newcommand{\R}{{\mathbb R}}

\newcommand{\cA}{{\cal A}}

\newcommand{\cO}{{\cal O}}
\newcommand{\cF}{{\cal F}}
\newcommand{\cG}{{\cal G}}

\newcommand{\cE}{{\cal E}}
\newcommand{\cH}{{\cal H}}

\newcommand{\cN}{{\cal N}}

\newcommand{\cL}{{\cal L}}

\newcommand{\cV}{{\cal V}}
\newcommand{\cW}{{\cal W}}

\newcommand{\g}{{\nabla}}

\newcommand{\pd}{\partial}
\newcommand{\intl}{\int\limits}

\newcommand{\Lto}{\widehat L_2(\Om)}
\newcommand{\Hto}{\widehat H^2_{0}(\Om)}
\newcommand{\hch}{\widehat{\cH}}
\newcommand{\wH}{\widehat{H}}

\newcommand{\supp}{{\operatorname{supp}}}

\newcommand{\di}{{\rm div\, }}
\newcommand{\wrt}{{with respect to }}

\newenvironment{declaration}[1]{\trivlist
\item[\hskip \labelsep{\bf #1 }]\ignorespaces}{\endtrivlist}
\newenvironment{proofof}[1]{\begin{declaration}{#1}}{\hfill
$\square$ \end{declaration}}
\newenvironment{proof}{\begin{proofof}{Proof.}}{\end{proofof}}

\begin{document}
\title{A global attractor  for   a  fluid--plate interaction model}
\author{Igor Chueshov\footnote{\small e-mail:
 chueshov@univer.kharkov.ua}
 ~~~ and ~~~
  Iryna Ryzhkova\footnote{\small e-mail:
 iryonok@gmail.com} \\
 \\Department of Mechanics and Mathematics, \\
 Kharkov National University, \\ Kharkov, 61077,  Ukraine\\  }

\maketitle

\begin{abstract}
We study asymptotic dynamics of a  coupled system consisting of linearized 3D
Navier--Stokes equations in a bounded domain  and a classical (nonlinear)
elastic plate equation for transversal displacement
on a flexible flat  part  of the boundary.
We show that this problem generates a semiflow on appropriate
phase space.
Our main result states the existence of a compact
finite-dimensional global attractor for this semiflow.
We do not assume any kind of mechanical damping in the plate component.
Thus our results means that dissipation of the energy in the fluid due to viscosity is sufficient to stabilize the system.
To achieve the result we first study the corresponding linearized model
and show that this linear model generates strongly continuous
exponentially stable semigroup.

\end{abstract}

\section{Introduction}

We consider a coupled (hybrid) system which describes
interaction of a homogeneous viscous incompressible fluid
which  occupies a domain $\cO$ bounded by
the (solid) walls of the container $S$ and a horizontal boundary $\Om$
on which a thin (nonlinear) elastic plate is placed.
The motion of the fluid is described
by  linearized 3D Navier--Stokes equations.
To describe  deformations of the plate
 we consider  a generalized plate model  which accounts only for
transversal displacements and
 covers
a general large deflection
Karman type model (see, e.g., \cite{lagnese,lag-2,lag-lions}
and also \cite{cl-book} and
the references therein).  However, our results can be also
applied in the cases of nonlinear Berger and Kirchhoff plates
(see the discussion in Section~\ref{sec:force}).
\par
This fluid--structure interaction
  model assumes that  large deflections of the plate produce
 small effect on the fluid. This  corresponds to
the case when the fluid fills the container which is large in comparison
with the size of the plate.
\par
We note that the mathematical studies of the problem
of fluid--structure interaction
in the case of viscous fluids and elastic plates/bodies
have a long history.
We refer to \cite{CDEG05,Ggob-jmfm08,Ggob-aa09,Ggob-mmas09,Kop98}
and the references therein for the case of plates/membranes,
to \cite{CS06} in the case of moving elastic bodies,
and to \cite{avalos-amo07,aval-tri07,aval-tri09,BGLT07,BGLT08,DGHL03}
in the case of elastic bodies with the fixed interface;
see also  the literature cited in these references.
\smallskip\par

Our mathematical model is formulated as  follows.
\par
 Let $\cO\subset \R^3$
 be a bounded domain  with a sufficiently smooth
 boundary $\partial\cO$. We assume that
$\partial\cO=\overline{\Omega}\cup \overline{S}$,
 where $\Om\cap S=\emptyset$ and
$$
\Om\subset\{ x=(x_1;x_2;0)\, :\,x'\equiv
(x_1;x_2)\in\R^2\}
$$ with the smooth contour $\Gamma=\partial\Om$
and $S$
 is a  surface which lies in the subspace $\R^3_- =\{ x_3\le 0\}$.
 The exterior normal on $\partial\cO$ is denoted
 by $n$. We have that $n=(0;0;1)$ on $\Om$.
 We consider the following {\em linear} Navier--Stokes equations in $\cO$
for the fluid velocity field $v=v(x,t)=(v^1(x,t);v^2(x,t);v^3(x,t))$
and for the pressure $p(x,t)$:
\begin{equation}\label{fl.1}
   v_t-\nu\Delta v+\nabla p=G_f(t)\quad {\rm in\quad} \cO
   \times(0,+\infty),
\end{equation}
   \begin{equation}\label{fl.2}
   \di v=0 \quad {\rm in}\quad \cO
   \times(0,+\infty),
  \end{equation}
where $\nu>0$ is the dynamical viscosity and $G_f(t)$ is a volume force
(which may depend on $t$).
   We supplement (\ref{fl.1}) and (\ref{fl.2}) with  the (non-slip)  boundary
   conditions imposed  on the velocity field $v=v(x,t)$:
\begin{equation}\label{fl.4}
v=0 ~~ {\rm on}~S;
\quad
v\equiv(v^1;v^2;v^3)=(0;0;u_t) ~~{\rm on} ~ \Om.
\end{equation}
Here $u=u(x,t)$ is the transversal displacement
of the plate occupying $\Om$ and satisfying
the following  equation  (see, e.g., \cite{Bol63,lagnese,lag-2,lag-lions}
and the references therein):
\begin{equation*}
u_{tt} + \De^2 u + \cF(u)=G_{pl}(t)- T_f(v)
~~{\rm in}~~ \Omega \times (0, \infty),
\end{equation*}
where $G_{pl}(t)$ is a given body force on the plate, $\cF(u)$ is a nonlinear feedback force which would be specified later  and $T_f(v)$ is a surface force exerted by the fluid on the plate,
$T_f(v)= (Tn\vert_\Om, n)_{\R^3}$, where $n$ is a outer unit normal to $\pd\cO$ at $\Om$  and $T=\{T_{ij}\}_{i,j=1}^3$ is the stress tensor of the fluid,
\[
T_{ij}\equiv T_{ij}(v)=\nu\left(v^i_{x_j}+v^j_{x_i}\right)-p\delta_{ij},
\quad i,j=1,2,3.
\]
Since  $n=(0;0;1)$ on $\Om$, we have that
 $T_f(v)=2\nu \pd_{x_3}v^3 -p$.
It also follows from \eqref{fl.2} and \eqref{fl.4} that
 $\pd_{x_3}v^3=0$  on $\Om$ and   thus we arrive at the equation
\begin{equation}\label{pl_eq}
u_{tt} + \De^2 u + \cF(u)=G_{pl}(t) +p|_\Om
~~{\rm in}~~ \Omega \times (0, \infty).
\end{equation}
We impose clamped boundary conditions on the plate
\begin{equation}
u|_{\pd\Om}=\left.\frac{\pd u}{\pd n} \right|_{\pd\Om}=0 \label{plBC}
\end{equation}
and supply \eqref{fl.1}--\eqref{plBC} with initial data of the form
\begin{equation}
 v(0)=v_0,\quad u(0)=u_0, \quad u_t(0)=u_1, \label{IC}
\end{equation}
We  note that  \eqref{fl.2} and \eqref{fl.4} imply the following
compatibility condition
\begin{equation}\label{Com-con}
\int_\Om u_t(x',t) dx'=0 \quad \mbox{for all}~~ t\ge 0.
\end{equation}
This condition fulfills when
\[
\int_\Om u(x',t) dx'=const \quad \mbox{for all}~~ t\ge 0,
\]
which can be interpreted as preservation of the volume of the fluid.
\par
We also note that a similar class of models was considered before
in \cite{Chu_2010,Ggob-jmfm08,Ggob-aa09,Ggob-mmas09}. The main difference
between \eqref{fl.1}--\eqref{IC}
and models in these  publications
is that the papers mentioned deal {\em only}  with
 longitudinal deformations of the plate
neglecting transversal deformations
(in contrast with the model  \eqref{fl.1}--\eqref{IC} which takes
into account the transversal deformations only).
This means that instead of (\ref{fl.4}) the following
boundary conditions are imposed on the velocity fluid field:
\begin{equation}\label{fl.4-old}
v=0 ~~ {\rm on}~S;
\quad
v\equiv(v^1;v^2;v^3)=(u^1_t; u_t^2;0) ~~{\rm on} ~ \Om,
\end{equation}
where $u=(u^1(x,t); u^2(x,t))$ is the
in-plane displacement vector of
the plate which solves the wave equation of the form
\begin{equation}\label{pl-1m}
                        u_{tt} -\Delta u- \g
\left[ {\rm div}\, u\right]
+\nu    (v^1_{x_3};v^2_{x_3})|_{x_3=0}  +
                        f( u)=0~~{\rm in}~~\Om;~~~u^i =0~~{\rm on}~~
  \Gamma.
\end{equation}
 This kind of models arises in
the study of blood flows in large arteries
(see the references in \cite{Ggob-jmfm08}).
The model \eqref{fl.1}, \eqref{fl.2}, (\ref{fl.4-old}), (\ref{pl-1m})
is simpler in several respects.
One of them is related to the fact the force exerted on the plate by the fluid
is more regular in the case (\ref{pl-1m}) and does not contains
the pressure  in an explicit form.
Moreover, the  model \eqref{fl.1}, \eqref{fl.2},
(\ref{fl.4-old}),  (\ref{pl-1m}) does not require any compatibility conditions like (\ref{Com-con}), because the volume of the fluid obviously preserves in the case of longitudinal deformations.
\medskip \par
In this paper our main point of interest is well-posedness and
long-time dynamics of
solutions to the coupled problem in
\eqref{fl.1}--\eqref{IC}
for the velocity $v$ and the displacement $u$.
First we consider  the linear version of this problem (i.e., the case
when $\cF(u)\equiv 0$). For this linear version
we prove well-posedness in the class of weak (energy) solutions and
establish some additional properties of solutions
which we need for treating the nonlinear problem.
In particular, we
show that in the homogeneous case
($G_f\equiv 0$, $G_{pl}\equiv 0$)  the linear version   generates
strongly continuous exponentially stable semigroup.
Then we consider a nonlinear version of this problem under
rather general hypotheses concerning nonlinearity.
These hypotheses cover the cases of von  Karman, Berger and Kirchhoff
plates. We show that  problem \eqref{fl.1}--\eqref{IC} generates
a dynamical system in an energy type space.
Our main result (see Theorem~\ref{th:attractor})
states that under some natural conditions concerning
 feedback forces system \eqref{fl.1}--\eqref{IC}
possesses a compact global attractor of finite fractal dimension.
To establish this results we rely on recently developed approach (see \cite{cl-jdde}, \cite{cl-mem} and
\cite[Chapters 7,8]{cl-book} and also the references therein)
which  involves stabilizability estimates
and notion of a quasi-stable system.
\par

\par
The paper is organized as follows.
 In Section~\ref{sec:pre} we introduce   notations, recall
some properties of Sobolev type spaces with non-integer
indexes   on bounded domains and collect some regularity
properties of (stationary)
Stokes problem which we use in the further considerations
(see Proposition~\ref{pr:stokes}).
Section \ref{sec:lin} is devoted to  a linear version of
the problem. Our main result in this section
is Theorem~\ref{lin_WP} on  well-posedness of weak solutions.
In Section~\ref{sec:nonlin} we deal with the nonlinear
problem   \eqref{fl.1}--\eqref{IC}.
First we prove  well-posedness result in
 Theorem ~\ref{th:wp} and then show that in the
case of autonomous forces the problem generates a
gradient dynamical system.
Our main result in this  section states existence of a
finite dimensional global
attractor and describes  some regularity properties
of the trajectories from the attractor.
The  argument is based on the quasi-stability
property established in Proposition~\ref{pr:qst}.

\section{Preliminaries}\label{sec:pre}
In this section we introduce Sobolev type spaces we need and
provide with some results concerning to  Stokes problem.

\subsection{Spaces and notations}
To introduce Sobolev spaces we follow approach presented
in \cite{Triebel78}.
\par
Let $D$ be a sufficiently smooth domain  and $s\in\R$.
We denote  by $H^s(D)$ the Sobolev space of order $s$
on a set $D$ which we define as restriction (in the sense of distributions)
 of the
space $H^s(\R^d)$ (introduced via Fourier transform).
We denote by $\|\cdot \|_{s,D}$ the norm in  $H^s(D)$
which we define by the relation
\[
\|u\|_{s,D}^2=\inf\left\{\|w\|_{s,\R^d}^2\, :\; w\in H^s(\R^d),~~ w=u ~~
\rm{on}~~D
    \right\}
\]
We also use the notation $\|\cdot \|_{D}=\|\cdot \|_{0,D}$
for the corresponding $L_2$-norm and, similarly, $(\cdot,\cdot)_D$ for the $L_2$
inner product.
We denote by $H^s_0(D)$ the closure of $C_0^\infty(D)$ in  $H^s(D)$
(\wrt  $\|\cdot \|_{s,D}$) and introduce the spaces
\[
H^s_*(D):=\left\{f\big|_D\, :\;  f\in H^s(\R^d), \;
{\rm supp }\, f\subset \overline{D}\right\},\quad s\in \R.
\]
Since the extension by zero of elements from $H^s_*(D)$ gives us an
element of $H^s(\R^d)$,
these spaces $H^s_*(D)$ can be treated not only as functional spaces defined
on $D$ (and contained in $H^s(D)$) but also as  (closed) subspaces
of $ H^s(\R^d)$. Below we need them to describe
boundary traces on $\Om\subset\partial \cO$.
We endow the classes $H^s_*(D)$ with the induced norms
 $\|f \|^*_{s,D}= \| f \|_{s,\R^d}$
for $f\in H^s_*(D)$. It is clear that
\[
\|f \|_{s,D}\le \|f \|^*_{s,D}, ~~ f\in H^s_*(D).
\]
It is known  (see \cite[Theorem 4.3.2/1]{Triebel78})
that $C_0^\infty(D)$ is dense in $H^s_*(D)$ and
\begin{align*}
& H^s_*(D)\subset H^s_0(D)\subset H^s(D),~~~ s\in\R;\\
& H^s_0(D) =  H^s(D),~~~ -\infty< s\le 1/2;\\
& H^s_*(D)= H^s_0(D),~~~ -1/2<  s<\infty,~~ s-1/2\not\in
\{ 0,1,2,\ldots\}.
\end{align*}
In particular, $H^s_*(D)= H^s_0(D)= H^s(D)$ for $|s|<1/2$. By
 \cite[Remark 4.3.2/2]{Triebel78} we also have that
 $H^s_*(D)\neq  H^s(D)$ for
 $|s|>1/2$. Note that  in the notations of \cite{LiMa_1968}
the space $H^{m+1/2}_*(D)$ is the same as $H^{m+1/2}_{00}(D)$ for every
 $m= 0,1,2,\ldots$ , and for $s=m+\si$ with $0<\si<1$ we have
\[
\|u \|^*_{s,D}=\left\{ \| u\|^2_{s,D}
+\sum_{|\al|=m}\int_D \,\frac{|D^\al u(x)|^2}{d(x,\pd D)^{2\si}}\, dx
\right\}^{1/2},
\]
where $d(x,\pd D)$ is the distance between $x$ and $\pd D$.
The norm  $\|\cdot \|^*_{s,D}$ is equivalent to
 $\|\cdot \|_{s,D}$ in the case when
$s>-1/2$ and  $s-1/2\not\in\{ 0,1,2,\ldots\}$,
but  not equivalent in general.
\par
Understanding adjoint spaces \wrt duality between
$C_0^\infty(D)$ and $[C_0^\infty(D)]'$
by Theorems 4.8.1 and 4.8.2 from \cite{Triebel78} we also have that
\begin{align*}
 [H^s_*(D)]'= H^{-s}(D),~ s\in\R, ~~~\mbox{and} ~~~
 [H^s(D)]' =  H_*^{-s} (D),~ s\in (-\infty,1/2).
\end{align*}
Below we also use the factor-spaces
$H^s(D)/\R$  with the naturally induced norm.
\medskip\par
To describe fluid velocity fields
we introduce the following spaces.
\par
Let $\mathscr{C}(\cO)$  be the class of
$C^\infty$ vector-valued solenoidal (i.e., divergence-free) functions
$ v=(v^1;v^2;v^3)$
on $\overline{\cO}$ which vanish in a neighborhood  of $S$
and such that $v^1=v^2=0$ on $\Om$.
We denote by $X$ the closure of $\sC(\cO)$ \wrt  the $L_2$-norm and
by $V$ the closure  \wrt the $H^1(\cO)$-norm. One
can see that
\[
X=\left\{ v=(v^1;v^2;v^3)\in [L_2(\cO)]^3\, :\; {\rm div}\, v=0;\;
\gamma_n v\equiv (v,n)=0~\mbox{on}~ S\right\}
\]
and
\[
V=\left\{
v=(v^1;v^2;v^3)\in [H^1(\cO)]^3\, \left| \begin{array}{l}
 {\rm div}\, v=0,\;
v=0~\mbox{on}~ S, \\ v^1=v^2=0~\mbox{on}~\Om \end{array} \right.
  \right\}.
\]
We equip   $X$ with $L_2$-type norm $\|\cdot\|_\cO$
and denote by $(\cdot,\cdot)_\cO$ the corresponding inner product.
The space $V$ is endowed  with the norm  $\|\cdot\|_V= \|\nabla\cdot\|_\cO$.
For some details concerning this type spaces we refer to \cite{temam-NS},
for instance.
\par
We also need the Sobolev spaces consisting of functions with zero average
on the domain $\Om$, namely
we consider the space
\[
\widehat{L}_2(\Om)=\left\{u\in L_2(\Om): \int_\Om u(x') dx' =0 \right\}
\]
and also  $\widehat H^s(\Om)=H^s(\Om)\cap\widehat L_2(\Om)$ for $s>0$
with the standard $H^s(\Om)$-norm.
The notations   $\widehat H^s_*(\Om)$ and $\widehat H^s_0(\Om)$
have a similar meaning.
\begin{remark}\label{re:hat-space}
{\rm Below we use $\widehat{H}^2_0(\Om)$ as a state space for the
displacement of the plate.
It is clear that $\widehat{H}^2_0(\Om)$ is a closed subspace of $H^2_0(\Om)$.
We denote by $\widehat{P}$   the projection
on  $\widehat{H}^2_0(\Om)$ in $H^2_0(\Om)$ which is orthogonal
with respect to the inner product $(\Delta\cdot, \Delta\cdot)_\Om$.
One can see that
$(I-\widehat{P})H^2_0(\Om)$
  consists of functions $u\in H^2_0(\Om)$ such that $\Delta^2u=const$ and
thus has dimension one.
}
\end{remark}

\subsection{Stokes problem}

In further considerations we need some regularity  properties
of the terms responsible for fluid--plate interaction.
To this end we consider
 the following Stokes problem
\begin{align}
  -\nu\Delta v+\nabla p= g, \quad
   \di v=0 \quad {\rm in}\quad \cO; \nonumber
\\
 v=0 ~~ {\rm on}~S;
\quad
v=(0;0;\psi) ~~{\rm on} ~ \Om,\label{stokes}
\end{align}
where $g\in [L^2(\cO)]^3$ and $\psi\in \Lto$ are given.
This type of boundary value problems for the Stokes equation
was studied by many authors  (see, e.g., \cite{lad-NSbook} and \cite{temam-NS}
and the references therein). We collect some properties of solutions
to \eqref{stokes} in the following assertion.
\begin{proposition}\label{pr:stokes}
With the reference to problem (\ref{stokes})
the following statements hold.
\begin{itemize}
\item [{ \bf (1)}]
Let $g\in [H^{-1+\sigma}(\cO)]^3$ and $\psi\in H^{1/2+\sigma}_*(\Om)$
be such that $\int_\Om\psi(x')dx'=0$.
Then
for every $0\le \sigma\le  1$
problem \eqref{stokes} has a unique solution
$\{v;p\}$ in $[H^{1+\sigma}(\cO)]^3\times[ H^{\sigma}(\cO)/\R]$ such that
\begin{equation}\label{stokes-bnd1}
\|v\|_{[H^{1+\sigma}(\cO)]^3}+\|p\|_{H^{\sigma}(\cO)/\R}
\le c_0\left\{\|g\|_{[H^{-1+\sigma}(\cO)]^3}+\|\psi\|_{H_*^{\sigma+1/2}(\Om)} \right\}.
\end{equation}
  \item [{ \bf (2)}]If $g=0$,  $\psi\in H_*^{-1/2+\si}(\Om)$,
$0\le \si\le 1$,  $\int_\Om\psi dx'=0$, then
  \begin{equation}\label{stokes-bnd2}
\|v\|_{[H^{\si}(\cO)]^3}+\|p\|_{H^{-1+\si}(\cO)/\R}
\le c_0\|\psi\|_{H_*^{-1/2+\si}(\Om)}.
\end{equation}
In particular, we can define a linear operator
$N_0 :\widehat{L}_2(\Om) \mapsto [H^{1/2}(\cO)]^3$
 by the formula
\begin{equation}\label{fl.n0}
N_0\psi=w ~~\mbox{iff}~~\left\{
\begin{array}{l}
 -\nu\Delta w+\nabla p=0, \quad
   \di w=0 \quad {\rm in}\quad \cO;
\\
 w=0 ~~ {\rm on}~S;
\quad
w=(0;0;\psi) ~~{\rm on} ~ \Om,
\end{array}\right.
\end{equation}
for $\psi\in \widehat{L}_{2}(\Om)$
($N_0\psi$ solves \eqref{stokes} with $g\equiv 0$).
It follows from (\ref{stokes-bnd1}) and (\ref{stokes-bnd2})
that
\[
N_0 :\, \widehat{H}^s_*(\Om)\mapsto [H^{1/2+s}(\cO)]^3\cap X
~~continuously~for~ -\frac12\le s\le \frac32.
\]
  \item [{ \bf (3)}] Let $g\in [H^{-1/2+\si}(\cO)]^3$ and $\psi\in \widehat{H}^{\si}_*(\Om)$, with $0<\si\le 1/2$.
Then we can define the trace of the pressure  $p$ on $\Om$,
which possesses the property  $p|_\Om\in H^{-1+\si}(\Om)/\R$ and
\begin{equation}\label{stokes-presure}
\|p\|_{H^{-1+\si}(\Om)/\R} \le c_0\left\{\|g\|_{[H^{-1/2+\si}(\cO)]^3}+\|\psi\|_{H^{\si}_*(\Om)} \right\}.
\end{equation}
\end{itemize}
\end{proposition}
\begin{proof}
Since the extension of elements from $H^\sigma_*(\Om)$ by zero
to the whole boundary $\partial\cO$  do not change
the smoothness Sobolev class, i.e., leads to elements from
$H^s(\partial\cO)$, we can use the regularity results
available for the Stokes problem with the Dirichlet type
boundary conditions imposed on the whole $\partial\cO$
(see, e.g., \cite{lad-NSbook,temam-NS} and also
the paper \cite{GSS2005} and the references therein).
This observation leads to the following arguments.
\par
{\bf 1.}
The existence and uniqueness of solutions along with the bound in
\eqref{stokes-bnd1} follow from Proposition 2.3 and Remark 2.6
on  Sobolev norm's interpolation in \cite[Chapter~1]{temam-NS}.

\par
{\bf 2.}
By Theorem 3\cite{GSS2005}  (applied for
the boundary data $\tilde\psi\in \wH^{-1/2}(\partial\cO)$
which is extension by zero outside $\Om$ of  the function $\psi\in \wH^{-1/2}_*(\Om)$)
 we have  (\ref{stokes-bnd2}) with $\si=0$.
Therefore interpolating with  (\ref{stokes-bnd1}) for  $s=0$
with $g\equiv 0$
we obtain (\ref{stokes-bnd2}) for all $0\le\si\le 1$.

\par
{\bf 3.}
We first represent $v$ in the form $v=\hat v+v^*$,
where $\hat v$ solves \eqref{stokes} with $\psi\equiv 0$
and $v^*$  satisfies \eqref{stokes}  with $g\equiv 0$.
Let $\hat p$ and $p^*$ be the corresponding representatives
of the pressure (which are identified with an element in a factor-space).
By the first statement
we have that $\hat p\in H^{1/2+\si}(\cO)$
and thus by the standard trace theorem there exists
$\hat p|_{\partial\cO}\in
H^{\si}(\partial\cO)$. This implies  that
$\hat p|_{\Om}\in H^{\si}(\Om)\subset  H^{-1+\si}(\Om) $ and
\begin{equation}\label{stokes-presure2}
\|\hat p\|_{H^{-1+\si}(\Om)/\R}\le c
\|\hat p\|_{H^{\si}(\Om)/\R}
\le c\|g\|_{[H^{-1/2+\si}(\cO)]^3}.
\end{equation}
In  the case $g\equiv 0$ the pressure $p^*$ is a harmonic function in $\cO$
which belongs $H^{-1/2+\si}(\cO)$.
This allows us to
 assign a meaning to $p^*|_\Om$ in $H^{-1+\si}(\Om)$.
Indeed, let $\phi\in C_0^\infty(\Om)$ and
 $\tilde\phi\in C^\infty_0(\partial\cO)$ be the extension
of $\phi$ by zero. Then by the trace theorem  there exists  a
 smooth function $w_\phi$ on $\cO$  such that
\[
w_\phi|_{\pd\cO}=0,\quad
\frac{\pd w_\phi}{\pd n}\Big\vert_{\pd\cO}=\tilde \phi, \quad \|w_\phi\|_{H^{5/2-\si}(\cO)}\le C\|\phi\|_{ H^{1-\si}_*(\Om)}.
\]
The application of  Green's formula yields
$(p^*,\Delta w_\phi)_\cO=
(p^*,\phi)_{\Om}$. Therefore
\[
|(\phi,p^*)_{\Om}|=|(p^*,\Delta w_\phi)_\cO|\le C \|p^*\|_{-1/2+\si,\cO}
\|\tilde\phi\|_{1-\si,\partial\cO}.
\]
Since
$\|\tilde \phi\|_{ 1-\si,\partial\cO}=\|\phi\|_{ H_*^{1-\si}(\Om)}$
and $C_0^\infty(\Om)$ is dense in  $H_*^{1-\si}(\Om)$,
we obtain
\begin{equation}\label{stokes-presure4}
\| p^*\|_{H^{-1+\si}(\Om)/\R}
\le c\|p^*\|_{H^{-1/2+\si}(\cO)/\R}\le c \|\psi\|_{H^{\si}_*(\Om)}.
\end{equation}
Thus relation (\ref{stokes-presure}) follows
from (\ref{stokes-presure2}) and (\ref{stokes-presure4}).
\end{proof}

\section{Linear problem}\label{sec:lin}
In this section we consider
 a linear version of \eqref{fl.1}--\eqref{IC}
which is  obtained from \eqref{fl.1}--\eqref{IC}
by replacing equation \eqref{pl_eq} with
its linear version. Thus we deal with  the following
problem
\begin{align}\label{fl.1-lin}
&   v_t-\nu\Delta v+\nabla p=G_f(t)~~\mbox{and}~~\di v=0~~{\rm in}~~
 \cO   \times(0,+\infty), \\[2mm]
& \label{fl.4-lin}
v=0 ~~ {\rm on}~S ~~\mbox{and}
~~
v\equiv(v^1;v^2;v^3)=(0;0;u_t) ~~{\rm on} ~~ \Om,
\\[2mm] &
u_{tt} + \De^2 u=G_{pl}(t)+p|_\Om  ~~{\rm on} ~~ \Om, \label{pl_eq-lin}
\\[2mm] &
u=\frac{\pd u}{\pd n} =0  ~~{\rm on} ~~ \pd\Om, \label{plBC-lin}
\end{align}
which we
supply  with the initial data of the form
\begin{equation}
v(0)=v_0,\quad u(0)=u_0, \quad u_t(0)=u_1.  \label{IC-lin}
\end{equation}
To define weak (variational) solutions
we need the following class $\cL_T$ of test functions $\phi$ on $\cO$:
\begin{equation*}
\cL_T=\left\{\phi \left|\begin{array}{l}
\phi\in L_2(0,T; \left[H^1(\cO)\right]^3),\; \phi_t\in L_2(0,T;  [L_2(\cO)]^3),  \\
{\rm div}\phi=0,\; \phi|_S=0,\; \phi|_\Om=(0;0;b),\; \phi(T)=0,  \\
b\in  L_2(0,T; \Hto),\; b_t\in  L_2(0,T; \Lto).
\end{array}\right.\right\}
\end{equation*}

\begin{definition}\label{lin_de:solution}
{\rm
A pair of functions $(v(t);u(t))$ is said to be  a weak solution to
the problem in
\eqref{fl.1-lin}--\eqref{IC-lin}  on a time interval $[0,T]$ if
\begin{itemize}
    \item $v\in L_\infty(0,T;X)\bigcap L_2(0,T; V)$;
\item $u \in L_\infty(0,T;H^2_0(\Om)), \; u_t \in L_\infty(0,T; \Lto)$
and   $u(0)=u_0$;
    \item for every $\phi\in \cL_T$  the following equality holds:
        \begin{multline}
            -\!\int_0^T\!\!(v,\phi_t)_\cO dt +\nu\!\int_0^T\!\!(\g v,\g\phi)_{\cO}  dt -\!\int_0^T\!\!(u_t,b_t)_\Om dt
           + \!\int_0^T\!\!(\De u, \De b)_\Om dt   \\
          =  \int_0^T(G_f(t), \phi)_{\cO} dt +\int_0^T(G_{pl}(t),b)_\Om dt + (v_0, \phi(0))_{\cO} + (u_1, b(0))_\Om; \label{weak_sol_def}
        \end{multline}
    \item the compatibility condition  $v(t)|_\Om=(0;0;u_t(t))$
holds for almost all $t$.
   \end{itemize}
}
\end{definition}
\begin{remark}\label{re:weak-sol}
{\rm {\bf (1)} It follows from the compatibility condition and
the standard trace theorem that
$u_t\in L_2(0,T; H^{1/2}_*(\Om))$ and
\begin{equation*}
||u_t(t)||_{H^{1/2}_*(\Om)} \le C||\g v(t)||_\cO
 \quad \mbox{for almost all } \; t\in [0,T].
\end{equation*}
{\bf (2)}   Taking in (\ref{weak_sol_def})
 $\phi(t)=\int_t^T\chi(\tau)d\tau\cdot \psi$, where $\chi$
is a smooth scalar function  and $\psi$ belongs to the space
\begin{equation}\label{space-W}
W=\left\{
\psi\in  V \left|  \;
  \psi|_\Om=(0;0;\beta),  \;
\beta\in \Hto \right. \right\},
\end{equation}
one see that the weak solution $(v(t);u(t))$  satisfies the relation
       \begin{multline}
           (v(t),\psi)_\cO
+(u_t(t),\beta)_\Om = (v_0, \psi)_{\cO} + (u_1, \beta)_\Om \\
-\int_0^t\left[
 \nu(\g v,\g\psi)_{\cO} +(\De u, \De \beta)_\Om
          -  (G_f, \psi)_{\cO}  -(G_{pl},\beta)_\Om\right] d\tau  \label{weak_sol_d2}
        \end{multline}
for almost all $t\in [0,T]$ and for all
$\psi=(\psi^1;\psi^2;\psi^3)\in W$, where $\beta=\psi^3\big\vert_\Om$.
}
\end{remark}
Below  as a phase space we  use
\begin{equation}\label{space-cH}
\cH=\left\{ (v_0;u_0;u_1)\in X\times H^2_0(\Om)\times\Lto :\; (v_0,n)\equiv
v_0^3 =u_1
~\mbox{on}~ \Om\right\}
\end{equation}
with the norm $\|(u_0;u_0;u_1)\|_\cH^2=\|v_0\|^2_{\cO}+\|\De u_0\|^2_{\Om}+\|u_1\|^2_{\Om}$. We also denote by $\hch$ a subspace in
$\cH$ of the form
\begin{equation}\label{space-cH-hat}
\hch= \left\{ (v_0;u_0;u_1)\in \cH :\; u_0\in \wH^2_0(\Om)\right\}.
\end{equation}
Our main result in this section is  the following well-posedness
theorem concerning the linear problem.
\begin{theorem} \label{lin_WP}
Assume that $U_0=(v_0;u_0;u_1)\in \cH$, $G_f(t)\in L_2(0,T; V')$
and  $G_{pl}(t)\in L_2(0,T; H^{-1/2}(\Om))$. Then
 for any interval $[0,T]$
there exists a unique weak solution $(v(t); u(t))$ to
\eqref{fl.1-lin}--\eqref{IC-lin}
 with the initial data $U_0$. This solution possesses the property
\begin{equation}\label{cont-ws}
U(t;U_0)\equiv U(t)\equiv (v(t); u(t); u_t(t))\in C(0,T; X\times H_0^2(\Om)\times \Lto),
\end{equation}
and satisfies  the energy balance equality
\begin{multline}\label{lin_energy}
\cE_0(v(t), u(t), u_t(t))+\nu \int_0^t ||\g v||^2_\cO d\tau=\cE_0(v_0, u_0, u_1) \\
+\int_0^t(G_f(\tau),  v)_\cO d\tau +\int_0^t (G_{pl}(\tau),u_\tau)_\Om d\tau
\end{multline}
 for every $t>0$, where the energy functional $\cE_0$ is defined
by the relation
\begin{equation}\label{lin_en-f}
\cE_0(v(t), u(t), u_t(t))=\frac12\left(\|v(t)\|^2_\cO+ \|u_t(t)\|^2_\Om+\| \De u(t)\|_\Om^2\right).
\end{equation}
Moreover, there exist positive constants $M$ and $\ga$ such that
for every initial data $U_0=(v_0;u_0;u_1)$ from  $\hch$ we have
\begin{equation}\label{exp-st-nhm}
\|U(t)\|^2_\cH\le M e^{-\ga t} \|U_0\|^2_\cH+M \int_0^t
e^{-\ga(t-\tau)}\left[ \|G_f(\tau)\|_{V'}^2+
\|G_{pl}(\tau)\|^2_{-1/2,\Om}\right] d\tau
\end{equation}
\end{theorem}

\begin{remark}\label{re:hat-space-2}
{\rm  Let $w_0\in (I-\widehat{P})H^2_0(\Om)$, where
the projector $\widehat{P}$ is defined in Remark~\ref{re:hat-space}.
Then one can see that the pair $\{v(t)\equiv 0, u(t)\equiv w_0\}$
solve  problem
\eqref{fl.1-lin}--\eqref{IC-lin}
 with the initial data $(0;w_0;0)$ and with $G_f\equiv 0$, $G_{pl}\equiv 0$.
The pressure $p$ is the constant determined from
its boundary value on $\Om$:
 $p|_\Om=\De^2w_0$ ($\De^2w_0$ is a constant due to Remark~\ref{re:hat-space}).
This observation gives us a relation
between solutions with initial data from $\cH$ and $\hch$, namely
we have that
\[
U(t; (v_0; u_0; u_1))= U(t; (v_0; \widehat{P} u_0; u_1))
+ (0;(I-\widehat{P})u_0; 0),
~~~ t>0,
\]
for any  $U_0=(v_0;u_0;u_1)\in \cH$. This  relation
means that $\hch$ is  invariant \wrt dynamics governed by \eqref{fl.1-lin}--\eqref{IC-lin}
and explains why an exponential decay estimate of the form
(\ref{exp-st-nhm}) cannot be true for {\em every} initial
data   $U_0=(v_0;u_0;u_1)$ from the space $\cH$.
}
\end{remark}
This remark allows us to derive from Theorem~\ref{lin_WP}
the following assertion.
\begin{corollary}\label{co:exp-stab}
Problem \eqref{fl.1-lin}--\eqref{IC-lin}
with $G_f\equiv 0$ and  $G_{pl}\equiv0$
generates a strongly continuous contraction semigroup $T_t$ on $\cH$ and on $\hch$ by the formula
    $T_t U_0=U(t)$, where $U(t)$ is a weak solution to \eqref{fl.1-lin}--\eqref{IC-lin}
with the initial data $U_0$. This semigroup $T_t$ is exponentially stable
on $\hch$, i.e.,
there exist positive constants $M$ and $\ga$ such that
\begin{equation*}
\|T_t U_0\|_\cH\le M e^{-\ga t} \|U_0\|_\cH~~~
for~ any ~~~ U_0=(v_0;u_0;u_1)\in \hch.
\end{equation*}
\end{corollary}

\begin{proof}
Strong continuity of $T_t$ follows from
(\ref{cont-ws}). This semigroup is contractive
and exponentially stable  due to  (\ref{lin_energy})
and (\ref{exp-st-nhm}) with   $G_f\equiv 0$ and $G_{pl}\equiv0$.
\end{proof}
We note that the generator of the semigroup $T_t$  defined via
solutions to problem
(\ref{fl.1-lin})--(\ref{IC-lin}) in the space $\hch$
has  a rather complicated structure, see Appendix~\ref{ap:generator}
in the end of the paper.
This is why we avoid in the argument below
calculations involving the explicit form of the generator.
\subsection*{Proof of Theorem~\ref{lin_WP}}
We use the compactness method and  split the argument
into  several steps.
\smallskip\par\noindent
{\em Step 1. Existence of an approximate solution.} For
the construction of Galerkin's approximations we use
an idea of  \cite{CDEG05} in a slightly modified form.
\par
Let $\{\psi_i\}_{i\in \N}$ be  the orthonormal basis in $\widetilde X
=\{v\in X : (v,n)\big\vert_\Om=0\}$
consisting of the eigenvectors of the Stokes problem:
\begin{gather*}
-\De \psi_i +\nabla p_i =\mu_i\psi_i  \quad \mbox{in} \; \cO,~~~
{\rm div}\psi_i=0, \quad \psi_i|_{\pd\cO}=0,
\end{gather*}
where  $0<\mu_1\le \mu_2\le \cdots$ are the corresponding eigenvalues.
 Denote by $\{\xi_i\}_{i\in\N}$ the basis in $\Hto$
which consists of eigenfunctions of the following problem
\[
(\De \xi_i,\De w)_\Om=\kappa_i (\xi_i, w)_\Om,~~~\forall\, w\in \wH^2_0(\Om),
\]
with the eigenvalues $0<\kappa_1\le\kappa_2\le\ldots$ and $||\xi_i||_{\Om}=1$.
Let  $\phi_i=N_0\xi_i$, where the operator $N_0$
is defined by (\ref{fl.n0}).
By Proposition~\ref{pr:stokes}  $\phi_i\in [H^2(\cO)]^3\cap V$.
As above one can also see that $\pd_{x_3}\phi^3_i=0$ on $\Om$.
 \par
We define an approximate solution as a pair of functions
\begin{equation}
v_{n,m}(t)=\sum_{i=1}^m \alpha_i(t)\psi_i +\sum_{j=1}^n \dot{\beta}_j(t)\phi_j, \quad u_n(t)=\sum_{j=1}^n\beta_j(t)\xi_j + (I-\widehat{P})u_0, \label{approx_sol}
\end{equation}
satisfying the relations
\begin{equation}
\dot{\alpha}_k(t) +\sum_{j=1}^n \ddot{\beta}_j(t)(\phi_j,\psi_k)_\cO+\nu\mu_k \alpha_k(t)+\nu\sum_{j=1}^n \dot{\beta}_j(t)(\g\phi_j,\g\psi_k)_\cO
=(G_f, \psi_k)_\cO  \label{psi_eq}
\end{equation}
for $k=1...m$, and
\begin{multline}
\sum_{i=1}^m \dot{\alpha}_i(t)(\psi_i, \phi_k)_\cO +\sum_{j=1}^n \ddot{\beta}_j(t)(\phi_j,\phi_k)_\cO+\ddot{\beta}_k(t)    \\ +
\nu\sum_{i=1}^m  \alpha_i(t)(\nabla\psi_i, \nabla\phi_k)_\cO +
\nu\sum_{j=1}^n \dot{\beta}_j(t)(\g\phi_j,\g\phi_k)_\cO
+\kappa_k \beta_k(t)  \\=
(G_f(t), \phi_k)_\cO +(G_{pl}(t), \xi_k)_\Om  \label{phi_eq}
\end{multline}
for $k=1,\dots,n$.
This system of ordinary differential equations is endowed
with the initial data
\[
v_{v,m}(0)=\Pi_m(v_0-N_0u_1)+N_0P_nu_1,
\]
\[
 u_n(0)=P_n\widehat{P}u_0 + (I-\widehat{P})u_0, \; \dot{u}_n(0)=P_n u_1,
\]
 where $\Pi_m$  is  the
 orthoprojector  on $Lin\{\psi_j : j=1,\ldots,m,\}$ in $\widetilde{X}$ and $P_n$
is orthoprojector on
$Lin\{\xi_i : i=1,\ldots,n\}$ in $\Lto$.
Since $\Pi_m$
and $P_n$ is are spectral projectors we have that
\begin{equation}\label{id-conv}
(v_{v,m}(0);u_n(0);  \dot{u}_n(0))\to  (v_0;u_0;u_1)~~
\mbox{strongly in $\cH$ as $m,n\to\infty$.}
\end{equation}

\par
We can rewrite system \eqref{psi_eq} and \eqref{phi_eq} as
\begin{equation*}
M\frac{d}{dt}\left(\begin{matrix}\alpha(t) \\ \dot{\beta}(t)\end{matrix}\right) +g(\alpha(t), \beta(t), \dot{\beta}(t))+G(t)=0
\end{equation*}
for some linear  function $g\,: \R^{m+2n}\mapsto \R^{m+n}$
and $G\in L_2(0,T;\R^{m+n})$, where
\begin{equation}\label{matr-M}
M=
\left[\begin{matrix}0&0\\0&  id\end{matrix}\right]+
\left[\begin{matrix} \{(\psi_i, \psi_j)_\cO\}_{j,k=1}^m & \{(\psi_l, \phi_k)_\cO\}_{l,k=1}^{m,n} \\ \{(\phi_k, \psi_l)_\cO\}_{l,k=1}^{n,m} &
\{ (\phi_i, \phi_j))_\cO\}_{j,k=1}^n\end{matrix}\right].
\end{equation}
The first matrix in (\ref{matr-M}) is nonnegative and
the second one is symmetric and   strictly  positive
(since the functions $\{\psi_i, \phi_j :i=1,\ldots,m, j=1,\ldots,n\}$ are linearly independent).
Therefore  system \eqref{psi_eq} and \eqref{phi_eq} has a unique  solution
on any time interval $[0,T]$.
\par
It follows from (\ref{approx_sol})  that
\[
v_{n,m}(t)=\sum_{i=1}^m \alpha_i(t)\psi_i + N_0[\pd_t u_n(t)],
\]
where $N_0$ is given by (\ref{fl.n0}). This implies
the following boundary compatibility condition
\begin{equation}\label{nm-comp}
v_{n,m}(t)=(0;0;\pd_t u_n(t))~~ \mbox{on}~~ \Om.
\end{equation}
{\em Step 2. Energy relation and a priori estimate
for an approximate solution.}
It follows from  \eqref{psi_eq} and \eqref{phi_eq}  that
the approximate solutions
satisfy the relation
       \begin{multline}
           (\dot v_{n,m}(t),\chi)_\cO
+(\ddot u_n(t), h)_\Om
+
 \nu(\g v_{n,m}(t),\g\chi)_{\cO} +(\De u_n(t), \De h)_\Om
\\
 =
            (G_f(t), \chi)_{\cO}  +(G_{pl}(t),h)_\Om \label{app-sol1}
        \end{multline}
for  $t\in [0,T]$ and for
every $\chi$ and $h$ of the form
\begin{equation*}
\chi(t)=\sum_{k=1}^{m'} \chi_k\psi_k +N_0h~~\mbox{with} ~~
h=\sum_{k=1}^{n'} h_k\xi_k,
\end{equation*}
where $m'\le m$ and $n'\le n$. Therefore taking $\chi=v_{n,m}$ we obtain
 the following energy balance relation for approximate solutions
\begin{align} \label{approx_est}
&\cE_0(v_{n,m}(t), u_n(t), \pd_t u_n(t)) + \nu \int_0^t\int_\cO |\g v_{n,m}|^2 dx d\tau  \\ &  =\cE_0(v_{n,m}(0), u_n(0), \pd_t u_n(0)) +\int_0^t(G_f,v_{n,m} )_\cO d\tau + \int_0^t(G_{pl},\pd_t u_n)_\Om d\tau.
\nonumber
\end{align}
This implies the following a priori estimate
\begin{equation}
\sup_{t\in [0,T]}\left\{ \|v_{n,m}(t)\|^2_\cO+ \|\Delta u_n(t)\|^2_\Om+
 \|\pd_t u_n(t))\|^2_\Om\right\} + \! \int_0^T\!\|\g v_{n,m}\|_\cO^2 d\tau\le C_T. \label{approx_est1}
\end{equation}
By the trace theorem from (\ref{nm-comp}) we also have that
\begin{equation} \label{h12_est}
\int_0^T \|\pd_t u_n(\tau))\|^2_{H_*^{1/2}(\Om)} d\tau =
\int_0^T\| v_{n,m}(\tau)\|_{1/2,\pd\cO}^2 d\tau\le C_T.
\end{equation}
{\em Step 3. Limit transition.}
By \eqref{approx_est1} the sequence $\{(v_{n,m}; u_n; \pd_t u_n)\}$
contains a subsequence  such that
\begin{align}
&(v_{n,m}; u_n; \pd_t u_n) \rightharpoonup (v; u; \pd_t u) \quad \ast\mbox{-weakly in } L_\infty(0,T;\cH);\label{uv-conv} \\
&u_n \rightarrow u \quad \mbox{strongly in } C(0,T; H^{2-\e}_0(\Om))
,~~ \forall\, \eps>0;
\label{u-strong} \\
&v_{n,m} \rightharpoonup v \quad \mbox{weakly in } L_2(0,T;V).   \label{v_conv}
\end{align}
To obtain (\ref{u-strong}) we use the Aubin-Dubinsky  theorem
(see, e.g., \cite[Corollary~4]{sim}).
By (\ref{h12_est}) we can also suppose that
\begin{align}
&\pd_t u_n \rightharpoonup \pd_t u \quad \mbox{weakly in } L_2(0,T; H^{1/2}_*(\Om));\label{ut-conv}\\
&v_{n,m} \rightharpoonup v \quad \mbox{weakly in } L_2(0,T; H^{1/2}(\partial\cO)).
   \label{v_conv-b}
\end{align}
One can  see from (\ref{app-sol1}) that
$(v_{n,m}; u_n; \pd_t u_n)(t)$ satisfies (\ref{weak_sol_def}) with
the test function $\phi$ of the form
\begin{equation}\label{phi-pq}
\phi=\phi_{p,q}=\sum_{i=1}^p \gamma_i(t)\psi_i +\sum_{j=1}^q\delta_j(t)\phi_j,
\end{equation}
where $p\le m$, $q\le n$ and $\gamma_i$, $\delta_j$
are scalar absolutely continuous functions on $[0,T]$
such that $\dot{\ga}_i,\dot{\delta}_j\in L_2(0,T)$ and    $\gamma_i(T)=\delta_j(T)=0$. Thus using (\ref{uv-conv})--
(\ref{v_conv}) we can pass to the limit and
show that $(v; u; \pd_t u)(t)$ satisfies (\ref{weak_sol_def})
with  $\phi=\phi_{p,q}$, where $p$ and $q$ are arbitrary.
By (\ref{id-conv}) and (\ref{u-strong}) we have $u(0)=u_0$
The compatibility condition \eqref{fl.4-lin} follows from  (\ref{nm-comp})
and (\ref{ut-conv}), (\ref{v_conv-b}).
\par
To conclude the proof of the existence of weak solutions
we only need to show that any function $\phi$ in $\cL_T$ can be approximate by
a sequence of functions of the form (\ref{phi-pq}). This can be done in the following way. We first approximate the corresponding boundary value of $b$
by a finite linear combination $h$ of $\xi_j$, then we  approximate the difference $\phi-N_0h$ (with $N_0$ define by (\ref{fl.n0}))
by finite linear combination of $\psi_k$.
\par
Thus the existence of weak solutions is proved.
One can also see from (\ref{approx_est}) and
from (\ref{uv-conv})--(\ref{v_conv})
that the constructed weak solution satisfies the
corresponding energy balance {\em inequality}.
\par
{\it Step 4. Uniqueness.} We use the same idea as in \cite{Lions_1969},
but with a slightly modified test function, see   (\ref{phi-s}).
\par
Let $U^j(t)=(v^j(t);u^j(t);u_t^j(t))$, $j=1,2$, be two different solutions to the problem in question with the same initial data. Then their difference  $U(t)=U^1(t)-U^2(t)=(v(t);u(t);u_t(t))$ satisfies the variational equality
\begin{equation}
-\int_0^T(v,\phi_t)_\cO  + \nu\int_0^T  (\g v, \g \phi)_\cO -\int_0^T  (u_t, \pd_t b)_\Om + \\
\int_0^T (\Delta u, \Delta b)_\Om =0  \label{diff_var}
\end{equation}
for all $\phi\in\cL_T$, $b=(\phi|_\Om)^3$.
Now for every $0<s<T$ we take
\begin{equation}\label{phi-s}
\phi(t)\equiv\phi^s(t)=\left\{
    \begin{aligned}
         &-\int_t^s d\tau \int_0^\tau d\zeta v(\zeta), && t<s, \\
         &0, && t\ge s,
    \end{aligned}
\right.
\end{equation}
as a test function. We denote
\[
\psi^s(t)=\pd_t\phi^s(t)= \int_0^t d\zeta v(\zeta)~~\mbox{and}~~
b^s(t)=(\phi^s(t)|_\Om)^3=-\int_t^s d\tau u(\tau).
\]
 Substituting $\phi^s$ into \eqref{diff_var}, we obtain
\begin{equation}
    -\int_0^s (v, \psi^s)_\cO +\nu\int_0^s  (\g\pd_t\psi^s, \g\phi^s)_\cO -\int_0^s  (u_t, u)_\Om +\int_0^s  (\pd_t\Delta b^s,\Delta b^s)_\Om =0. \label{diff_var_1}
\end{equation}
Integrating by parts the second term in \eqref{diff_var_1}
and using the relations $\psi^s(0)=0$ and $\phi^s(s)=0$, we have
\[
\int_0^s  (\g\pd_t\psi^s, \g\phi^s)_\cO=
(\g\phi_s, \g\psi_s)_\cO\Big|_0^s - \int_0^s  (\g\psi^s,  \g\psi^s)_\cO
=-\int_0^s  \|\g\psi^s\|^2_\cO.
\]
Therefore \eqref{diff_var_1} yields
\begin{equation*}
||\psi^s(s)||^2_\cO+2\nu\int_0^s  \|\g\psi^s\|^2_\cO dt+ ||u(s)||^2_\Om + ||\Delta b^s(0)||^2_\Om=0
\end{equation*}
for almost all $0\le s\le T$. Therefore $v(s)=0$ and $u(s)=0$ for almost all $0\le s\le T$. Thus the uniqueness  is proved.
\par
{\it Step 5. Continuity with respect to $t$ and the energy equality.}
First we note that any weak solution $(v(t);u(t);u_t(t))$ is weakly continuous in
$X\times H^2_0(\Om) \times \Lto$. Indeed, it follows from (\ref{weak_sol_d2}) that
that any weak solution $(v(t);u(t))$  satisfies the relation
       \[
           (v(t),\psi)_\cO
 = (v_0, \psi(0))_{\cO}
+\int_0^t\left[-
 \nu(\g v,\g\psi)_{\cO}
          +  (G_f(\tau), \psi)_{\cO}  \right] d\tau
        \]
for almost all $t\in [0,T]$ and for all
$\psi\in \widetilde{V}=\{ v\in V :\, v|_\Om=0\}\subset W$.
This implies that $v(t)$ is weakly continuous in $\widetilde{V}'$.
Since $X\subset\widetilde{V}'$, we can apply the Lions lemma
(see \cite[Lemma 8.1]{LiMa_1968}) and conclude that $v(t)$ is
weakly continuous in $X$. The same lemma gives us weak continuity
 of $u(t)$ in $H^2_0(\Om)$. Now using  (\ref{weak_sol_d2}) again
we conclude that $(u_t(t),\beta)_\Om$ is continuous for $\beta\in \wH^2_0(\Om)$.
The density argument yields  weak continuity
 of $u_t(t)$ in $\Lto$.
\par
To prove the energy equality, we follow the scheme
of~\cite[Ch.1]{Lions_1969},
see also  \cite[Ch.3]{LiMa_1968}.
We  first note that due  to Remark~\ref{re:hat-space-2} it is
sufficient to consider the case when $U_0=(v_0;u_0;u_1)\in\hat\cH$.
Then
for every fixed $0<s<t<T$ we introduce
 a piecewise-linear continuous function $\theta_n(\tau)$
 on $\R$  such that  $\theta_n(\tau)=1$ for $s\le\tau\le t$ and   $\theta_n(\tau)=0$ when $\tau<s-1/n$ or $\tau>t+1/n$.
Let $\rho_k\in  C_0^\infty(\R)$ be an even function such that
$\supp\, \rho_k\subset [-k^{-1}, k^{-1}]$ and
$\int_\R\rho_k(s)ds=1$. Now for
 $k$ and $n$ large enough we consider
 the function $\phi(\tau)=\theta_n((\theta_n v) \ast \rho_k \ast \rho_k)$, where $v$ is a weak solution to \eqref{fl.1-lin}--\eqref{IC-lin},  as a test function in variational equality \eqref{weak_sol_def}.
Substituting this $\phi$ into \eqref{weak_sol_def} and passing to the limit
when $k\to\infty$ we obtain that
       \begin{multline}\label{ener-nn}
            -\!\int_0^T \tht_n\dot\tht_n\left[ \|v\|^2_\cO +
\|u_\tau\|^2_\Om+\|\De u\|^2_\Om \right] d\tau
+\nu\!\int_0^T \tht_n^2\|\g v\|^2_{\cO}  d\tau    \\
          =  \int_0^T \tht_n^2\left[ (G_f(\tau), v)_{\cO}
 +(G_{pl}(\tau),u_k)_\Om\right] d\tau
        \end{multline}
As in ~\cite[Ch.~1]{Lions_1969} one can see that for every
function $h\in L_1(0,T)$
\[
\lim_{n\to\infty}\int_0^T \tht_n(\tau)\dot\tht_n(\tau) h(\tau) d\tau=
-\frac12\left[ h(t)-h(s)\right]
\]
 for almost all $s$ and $t$. Therefore after the limit transition in
(\ref{ener-nn}) we obtain  energy relation (\ref{lin_energy})
 valid for almost all  $s$ and $t$. Now using weak continuity of
the solution $(v(t);u(t))$ and the energy  inequality
(which is valid for $s=0$ and for every $t$)
we can establish  the energy equality.
As in~\cite[Ch.~3]{LiMa_1968} this also implies
 strong continuity of weak solutions with respect to $t$.
\par
{\it Step 6. Exponential stability}. To prove the exponential stability
estimate   in (\ref{exp-st-nhm}), we construct a Lyapunov function
using an idea from \cite{Chu_2010}. Let
\begin{equation*}
V(v_0,u_0,u_1)=\cE_0(v_0,u_0,u_1)+\e\Psi(v_0,u_0,u_1),
\end{equation*}
where
$\Psi(v_0,u_0,u_1)=(u_0,u_1)_\Om +(v_0,N_0u_0)_\cO$
with $N_0$ defined by \eqref{fl.n0}, and $\e>0$ is
a small parameter which will be chosen later.
We consider these functionals on approximate solutions $(v_{n,m};u_n)$
for which $\widehat{P}u_0=u_0$ and thus  $\widehat{P}u_n(t)=u_n(t)$ for all $t>0$. This allow us to substitute in (\ref{app-sol1}) $N_0u_n$ instead of
$\chi$ and obtain that
  \begin{multline}
\frac{d}{dt}\Psi_{n,m}(t)\equiv
\frac{d}{dt}\Psi(v_{n,m}(t),u_{n}(t),
\pd_t u_{n}(t))=\|\pd_t u_{n}\|^2_\Om + (v_{n,m}, N_0\pd_tu_n)_{\cO}
\\
-
 \nu(\g v_{n,m},\g N_0u_n)_{\cO} -\|\De u_n\|^2_\Om
           + (G_f, N_0u_n)_{\cO}  +(G_{pl},u_n)_\Om. \label{app-psi}
        \end{multline}
By Proposition~\ref{pr:stokes}, using the compatibility
condition in (\ref{nm-comp}) and the trace theorem  we have that
\[
| (v_{n,m}, N_0\pd_t u_n)_{\cO}|\le C \|v_{n,m}\|_\cO \|\pd_tu_n\|_{\Om}
\le C  \|\g v_{n,m}\|^2_\cO.
\]
Similarly,
\[
  | (\g v_{n,m},\g N_0u_n)_{\cO}|\le \eta \|\De u_n\|^2+ C_\eta
   \|\g v_{n,m}\|^2_\cO,~~~\forall\, \eta>0,
\]
and also
\[
| (G_f, N_0u_n)_{\cO}  +(G_{pl},u_n)_\Om|\le  \eta \|\De u_n\|^2+C_\eta \left[\|G_f\|_{V'}^2 +\|G_{pl}\|^2_{-1/2,\Om}\right].
\]
Therefore it follows from (\ref{app-psi}) that
\[
\frac{d}{dt}\Psi_{n,m}(t)\le -\frac12 \|\De u_n\|^2+ C  \|\g v_{n,m}\|^2_\cO
+C\left[\|G_f\|_{V'}^2 +\|G_{pl}\|^2_{-1/2,\Om}\right].
\]
Using the energy relation in (\ref{approx_est}) we also have that
\begin{align*}
&\frac{d}{dt}\cE_0(v_{n,m}(t), u_n(t), \pd_t u_n(t)) \le -\frac{\nu}2
\|\g v_{n,m}\|^2_{\cO}
+C_\nu \left[\|G_f\|_{V'}^2 +\|G_{pl}\|^2_{-1/2,\Om}\right].
\end{align*}
Therefore the function
$V_{n,m}(t)\equiv V(v_{n,m}(t),u_{n}(t),
\pd_t u_{n}(t))$  satisfies the relations
\[
a_0 \cE_0(v_{n,m}(t),u_{n}(t),
\pd_t u_{n}(t))\le
V_{n,m}(t)\le a_1 \cE_0(v_{n,m}(t),u_{n}(t),
\pd_t u_{n}(t))
\]
for sufficiently small $\eps>0$  and
\[
\frac d{dt} V_{n,m}(t)+a_2 V_{n,m}(t)\le a_3  \left[\|G_f\|_{V'}^2 +\|G_{pl}\|^2_{-1/2,\Om}\right]
\]
with positive constants $a_i$. This implies relation (\ref{exp-st-nhm})
for approximate solutions. The limit transition yields  (\ref{exp-st-nhm})
for every weak solutions.
\par
This completes the proof of Theorem~\ref{lin_WP}.

\section{Nonlinear problem}\label{sec:nonlin}
In this section we deal with problem (\ref{fl.1})--(\ref{IC})
with a nonlinear feedback force.
Fist we  describe hypotheses concerning this force.
Then we prove well-posedness (see Theorem~\ref{th:wp})
and construct the corresponding semiflow.
Our main result (see Theorem~\ref{th:attractor}) states
the existence of finite-dimensional attractor.
\subsection{Structure of feedback force}\label{sec:force}
We impose the following hypotheses concerning the nonlinear feedback force $\cF(u)$ in the plate equation (\ref{pl_eq}).

\begin{assumption}\label{A:force}
\begin{itemize}
 \item[{\bf (F1)}]
There exists $\e>0$ such that $\cF(u)$ is locally Lipschitz from $H^{2-\e}_0(\Om)$ into $H^{-1/2}(\Om)$\footnote{
We recall that according our definitions $ H^{-1/2}(\Om) =[H^{1/2}_*(\Om)]'
 \varsupsetneqq [H^{1/2}_0(\Om)]'$.}
 in the sense that
\begin{equation}\label{f-lip}
\| \cF(u_1)-\cF_2(u_2)\|_{-1/2,\Om}\le C_R \| u_1-u_2\|_{2-\eps,\Om}
\end{equation}
for any $u_i\in H^{2}_0(\Om)$ such that $\| u_i\|_{2,\Om}\le R$.
  \item[{\bf (F2)}]
There exists a $C^1$-functional  $\Pi(u)$  on $H^2_0(\Om)$
such that $\cF(u)=\Pi'(u)$, where $\Pi'$ denotes the
 Fr\'echet derivative of  $\Pi$.
   \item[{\bf (F3)}]
The plate force potential $\Pi$ is bounded on bounded sets from $H^2_0(\Om)$
 and  there exist $\eta<1/2$
 and $C\ge 0$ such that
\begin{equation}\label{8.1.1c1}
 \eta \|\De u\|_\Om^2 +\Pi(u)+C \ge 0\;,\quad \forall\, u\in H^2_0(\Om).
\end{equation}
\end{itemize}
\end{assumption}
 The nonlinear feedback (elastic) force $\cF(u)$ may have one
of the following forms (which represent different plate models):
\medskip\par\noindent
 {\bf Kirchhoff model}: $\cF(u)$ is the Nemytskii operator
\begin{equation*}
u\mapsto - \kappa\cdot {\rm div}\left\{|\nabla u|^q\nabla u
-\mu |\nabla u|^r\nabla u
\right\}+ f(u)-h(x),
\end{equation*}
 where  $\kappa\ge 0$, $q>r\ge 0$ and $\mu\in \R$ are parameters, $h\in L_2(\Om)$, and
 \begin{equation}
    \label{phi_condition-r}
f\in {\rm Lip_{loc}}(\R)~~ \mbox{ satisfies }  ~~  \underset{|s|\to\infty}{\liminf}\, f(s)s^{-1}>-\la_1,
\end{equation}
 where $\la_1$ is the first eigenvalue of the
biharmonic operator with the Dirichlet boundary conditions.
In this case the relation in  (\ref{f-lip})  follows from the considerations given in \cite[Sect.5]{ChuKol2}. We also have that
\begin{align*}
\Pi(u) =& \int_\Om F(u(x))dx +\frac{\kappa}{q+2}\int_\Om |\nabla u(x)|^{q+2}dx
\\ &
-\frac{\kappa\mu}{r+2}\int_\Om |\nabla u(x)|^{r+2}dx
 -\int_\Om u(x) h(x)dx,
\end{align*}
where $F(s)=\int_0^sf(\xi)d\xi$ is the antiderivative of $f$.
Due to the second relation in (\ref{phi_condition-r}) we obviously
have (\ref{8.1.1c1}).
\medskip\par \noindent
 {\bf Von Karman model:} This model is
 well known in nonlinear elasticity and
constitute a basic model describing nonlinear oscillations of a
plate accounting for  large deflections, see \cite{ Lions_1969,cl-book}  and the
references therein. The force $\cF$ has the form
$\cF(u)=-[u, v(u)+F_0]-h(x)$, where
 $F_0\in H^4(\Om) $ and $h\in L_2(\Om)$ are given functions,
the von Karman bracket $[u,v]$  is given by
\begin{equation*}
[u,v] = \partial ^{2}_{x_{1}} u\cdot \partial ^{2}_{x_{2}} v +
\partial ^{2}_{x_{2}} u\cdot \partial ^{2}_{x_{1}} v -
2\cdot \partial ^{2}_{x_{1}x_{2}} u\cdot \partial ^{2}_{x_{1}x_{2}}
v ,
\end{equation*} and
the Airy stress function $v(u) $ solves the following  elliptic
problem
\begin{equation*}
\Delta^2 v(u)+[u,u] =0 ~~{\rm in}~~  \Omega,\quad \frac{\pd v(u)}{\pd
n} = v(u) =0 ~~{\rm on}~~  \pd\Om.
\end{equation*}
It is known (see, e.g., Corollary 1.4.5  in \cite{cl-book})
that
\[
\| [u_1,v(u_1)]-  [u_2,v(u_2)]\|_{-\eta,\Om}\le
C(\|u_1\|_{2,\Om}^2+\|u_2\|^2_{2,\Om})\| u_1-  u_2\|_{2-\eta,\Om}
\]
for every $\eta\in [0,1]$,
which implies (\ref{f-lip}). The potential  $\Pi$ has the form
\[
\Pi(u)=\frac14\int_\Om \left[ |v(u)|^2 -2([u,F_0]-2 h) u\right] dx
\]
and possesses the properties listed in Assumption~\ref{A:force},
see, e.g., \cite[Chapter 4]{cl-book} for details.
\medskip\par\noindent
 {\bf Berger Model:} In this case the feedback force has the form
  $$
  \cF(u)=- \left[ \kappa \int_\Om |\nabla u|^2 dx-\Gamma \right] \De u -h(x),
$$
where $\kappa>0$ and $\Gamma\in\R$ are parameters,  $h\in L_2(\Om)$.
One can see Assumption~\ref{A:force} is satisfied,
 for some details and  references see, e.g.,
\cite[Chapter 4]{Chueshov} and \cite[Chapter 7]{cl-mem}.

\subsection{Well-Possedness}\label{sec:wp}

\begin{definition}\label{de:solution}
{\rm
A pair of functions $(v(t);u(t))$ is said to be  a weak solution to \eqref{fl.1}--\eqref{IC}  on a time interval $[0,T]$ if
\begin{itemize}
    \item $v\in L_\infty(0,T;X)\bigcap L_2(0,T; V)$;
\item $u \in L_\infty(0,T;H^2_0(\Om)), \; u_t \in L_\infty(0,T; \Lto)$,  $u(0)=u_0$;
    \item the  equality
in (\ref{weak_sol_def}) holds with $G_{pl}(t):= -\cF(u(t))+  G_{pl}(t)$;
    \item the compatibility condition  $v(t)|_\Om=(0;0;u_t(t))$ holds for almost all $t$.
\end{itemize}
}
\end{definition}

\begin{theorem}\label{th:wp}
Assume that $U_0=(v_0;u_0;u_1)\in \cH$, $G_f(t)\in L_2(0,T; V')$
and  $G_{pl}(t)\in L_2(0,T; H^{-1/2}(\Om))$. Then
 for any interval $[0,T]$
there exists a unique weak solution $(v(t); u(t))$ to
\eqref{fl.1}--\eqref{IC}
 with the initial data $U_0$. This solution possesses the property
\begin{equation}\label{cont-ws-n}
U(t)\equiv (v(t); u(t); u_t(t))\in C(0,T; \cH),
\end{equation}
where $\cH$ is given by (\ref{space-cH}),
and satisfies  the energy balance equality
\begin{multline}\label{energy}
\cE(v(t), u(t), u_t(t))+\nu \int_0^t ||\g v||^2_\cO d\tau=\cE(v_0, u_0, u_1) \\
+\int_0^t(G_f(\tau),  v)_\cO d\tau +\int_0^t (G_{pl}(\tau), u_\tau)_\Om d\tau
\end{multline}
 for every $t>0$, where the energy functional $\cE$ is defined
by the relation
\begin{equation*}
\cE(v, u, u_t)=\frac12\|v\|^2_\cO+E(u, u_t)
\end{equation*}
with the plate energy  $ E(u,u_t)$ given by
\[
E(u,u_t)=\frac12\left( \|u_t\|^2_\Om+\| \De u\|_\Om^2\right)+ \intl_\Om\Pi(u(x))dx.
\]
Moreover, there exists a constant $a_{R,T}>0$ such that
  for any couple of weak solutions
  $U(t)=(v(t); u(t); u_t(t))$ and $\hat U(t)=(\hat v(t); \hat u(t); \hat u_t(t))$
with the initial data  possessing the property  $\|U_0\|_\cH, \|\hat U_0\|_\cH\le R$
we  have
\begin{equation}\label{mild-dif}
\|U(t)-\hat U(t)\|^2_\cH+\intl_0^t\|\nabla( v-\hat v)\|_\cO^2 d\tau \le a_{R,T}
\|U_0-\hat U_0\|^2_\cH, ~~t\in [0,T].
\end{equation}
The spatial average of $u(t)$  is preserved. In particular, if $U_0\in\hch$,
then $U(t)\in\hch$ for every $t>0$.
We recall that $\hch$ is defined by (\ref{space-cH-hat}).
\end{theorem}

\begin{proof}
The proof of the local existence of an approximate solution is almost the same, as in the linear case (see Theorem \ref{lin_WP}).
We use approximate solutions of the same structure as in
(\ref{approx_sol}) which satisfy (\ref{psi_eq}), (\ref{id-conv}) and also
(\ref{phi_eq})  with $-\cF(u_n(t))+  G_{pl}(t)$
instead of $G_{pl}(t)$.
Then using  the standard argument we establish
the energy relation in (\ref{energy}) for these approximate solutions.
Now
the positivity type estimate in (\ref{8.1.1c1}) allow us
to obtain  the same a priori estimates as
in (\ref{approx_est1}) and (\ref{h12_est}). Therefore
we can prove the global existence of approximate solutions and
 establish  the existence
of a weak solution $U(t)=(v(t); u(t); u_t(t))$ by the same
argument as in the linear case.
To  make  limit transition in the nonlinear term we
use (\ref{f-lip}).
\par
Now we can consider the  pair $(v(t); u(t))$ as a solution
to linear problem with $G_{pl}(t):= -\cF(u(t))+  G_{pl}(t)$.
This allow us to obtain (\ref{cont-ws-n}) and also derive
energy balance relation  (\ref{energy}) from (\ref{lin_energy})  using
the potential structure of the force $\cF$:  $\cF(u)=\Pi'(u)$.
\par
Since the difference of two weak solution  can be treated as
a solution to the linear problem with $G_f\equiv 0$
and  $G_{pl}(t):= \cF(\hat u(t)) -\cF( u(t))$,
 we can obtain (\ref{mild-dif}) from the energy
equality  (\ref{lin_energy}).  The uniqueness
follows from  (\ref{mild-dif}).
\par
Preservation of the spatial average of $u(t)$
follows from the same property for approximate solutions.
\end{proof}
\begin{remark}\label{re:autn}
{\rm In the autonomous case we can suggest another form of
energy relation (\ref{energy}).
Let $G_{pl}(t)\equiv 0$  and
 $G_f(t)\equiv G_0\in V'$ be independent of $t$. Suppose that
a pair $(v_*;p_*)\in V\times L^2(\cO)$ solve  problem (\ref{stokes}) with $g\equiv G_0$   and $\psi\equiv 0$, i.e.,
\begin{equation}\label{v-stat}
  -\nu\Delta v_*+\nabla p_*= G_0, ~~
   \di v_*=0 ~~ {\rm in} ~ \cO;
\quad v_*=0 ~~ {\rm on}~\pd\cO.
\end{equation}
Then the following
form of the energy  balance equation is valid:
\begin{equation}\label{en-shift}
\cE_*(v(t), u(t), u_t(t))+\nu \int_0^t ||\g (v-v_*)||^2_\cO d\tau=
\cE_*(v_0, u_0, u_1),
\end{equation}
where
\begin{equation*}
\cE_*(v, u, u_t)=\frac12\|v-v_*\|^2_\cO+E_*(u, u_t)
\end{equation*}
with $ E_*(u,u_t)$ given by
\[
E_*(u,u_t)=\frac12\left( \|u_t\|^2_\Om+\| \De u\|_\Om^2\right)+ \intl_\Om\Pi(u(x))dx -(p_*,u)_\Om.
\]
Indeed, it follows from (\ref{v-stat}) that
\[
(G_0,v(t))_\cO=\nu (\g v_*, \g v(t))_\cO+\frac{d}{dt}(p_*,u(t))_\Om.
\]
Substituting $\psi=v_*$ in (\ref{weak_sol_d2}) we also have that
\[
\frac{d}{dt}(v(t), v_*)_\cO+
\nu (\g v_*, \g v(t)) =(G_0, v_*)_\cO =\nu \|\g v_*\|^2_\cO.
\]
Therefore
\[
(G_0,v(t))_\cO=\frac{d}{dt}\left[ (v(t), v_*)_\cO+(p_*,u(t))_\Om\right]+
2\nu (\g v_*, \g v(t)) -\nu \|\g v_*\|^2_\cO.
\]
This and also the energy relation in (\ref{energy}) imply
(\ref{en-shift}).
}
\end{remark}
This remark allows us the derive from  Theorem~\ref{th:wp} the following
assertion.
\begin{corollary}\label{co:generation}
Let  $G_f(t)\equiv G_0\in V'$ be independent of $t$
 and $G_{pl}(t)\equiv 0$.
Then   problem \eqref{fl.1}--\eqref{IC} generates dynamical systems
 $(S_t, \cH)$ and $(S_t, \widehat{\cH})$ with the evolution operator
defined  by the formula $S_t U_0=(v(t);u(t);u_t(t))$,
where $(v;u)$ is a weak  solution to \eqref{fl.1}--\eqref{IC}
 with the initial data $U_0=(v_0; u_0; u_1)$.
These systems are gradient with the full energy $\cE_*(v_0, u_0, u_1)$
as a Lyapunov function. This means that
(a) $U\mapsto \cE_*(U)$
is continuous on $\cH$,
(b)
 $\cE_*(S_t U_0)$ is not increasing in $t$, and
(c)
if $\cE_*(S_t U_0)=\cE_*(U_0)$ for some $t>0$, then
 $U_0$ is a stationary point of $S_t$
(i.e., $S_t U_0=U_0$ for all $t\ge 0$).
Moreover,
the set $\cE_R=\{U_0: \cE_*(U_0)\le R\}$ is a bounded
closed forward invariant set for every $R>0$.
\end{corollary}
\begin{proof}
We need only to check the properties of the functional $\cE_*$.
\par
It is clear from Assumption~\ref{A:force}(F2) that $\cE_*$
is continuous on $\cH$.
\par
By (\ref{en-shift}) We have that
 $\cE_*(S_t U_0)\le \cE_*(S_\tau U_0)$ for $t\ge \tau\ge 0$.
This gives the monotonicity of  $t\mapsto \cE_*(S_t U_0)$
and the invariance of $\cE_R$.
\par
If  $\cE_*(S_{t_0} U_0)= \cE_*(U_0)$ for  some $t_0>0$, then
(\ref{en-shift}) implies that $v(t)=v_*$ for all $t\in [0,t_0]$
and thus $u_t(t)=v_*^3\big|_\Om=0$. Hence $u(t)\equiv u$ for
some $u\in H^2_0(\Om)$ and $U_0=(v_*; u; 0)$ is a stationary
point for $S_t$.
\end{proof}
Below we describe the set of stationary point
of the evolution semigroup $S_t$ with more details.

\subsection{Stationary solutions}
As above we assume that $G_{pl}\equiv 0$ and
$G_f(t)\equiv G_0\in V'$ is independent of $t$.
Let $\widetilde{V}=\{ u\in V\,: \, v\big|_{\pd\cO}=0\}$.
It follows from Definition~\ref{de:solution} that a stationary
(time-independent) solution  is a pair $(v;u)$
 from $\widetilde{V}\times H_0^2(\Om)$
satisfying the relation
       \begin{equation}
 \nu(\g v,\g\psi)_{\cO} +(\De u, \De \beta)_\Om
          -  (G_0, \psi)_{\cO}  +(\cF(u),\beta)_\Om=0
\label{weak_sol_stat}
        \end{equation}
for any $\psi\in W$ with $\psi^3\big\vert_\Om=\beta$,
where $W$ is given by (\ref{space-W}).
Using  (\ref{en-shift}) we have that $v=v_*$, where $v_*$ solves
(\ref{v-stat}).
One can also  see  $(\g v,\g N_0\beta)_\cO=0$ for any $v\in V_0$
and $\beta\in \wH^2_0(\Om)$, where $N_0$ is defined in~(\ref{fl.n0}).
Therefore from (\ref{weak_sol_stat}) with $\psi=N_0\beta$ we have the
following variational problem for $u\in H^2_0(\Om)$:
      \begin{equation}
 (\De u, \De \beta)_\Om
        +(\cF(u)-N_0^*G_0,\beta)_\Om=0,~~~\forall\, \beta\in\wH^2_0(\Om).
\label{plate_sol_stat}
        \end{equation}
The following calculation performed first on smooth functions  gives
us
\begin{multline*}
(G_0,N_0\beta)_\cO= ( -\nu\Delta v_*+\nabla p_*,N_0\beta)_\cO
\\
=
 (v_*, -\nu\Delta N_0\beta)_\cO + (p_*,\beta)_\Om
   = (v_*, -\g p_\beta)_\cO + (p_*,\beta)_\Om = (p_*,\beta)_\Om.
\end{multline*}
Since the pressure $p_*$ in (\ref{v-stat}) is defined up to a constant,
we can suppose that $p_*=N_0^*G_0$. By Proposition~\ref{pr:stokes}
$N_0^*\, :\ V'\mapsto [\wH^{1/2}_*(\Om)]'$. This provides us with the
regularity of the pressure impact on the plate.
\par
One can see that a function $u\in H^2_0(\Om)$ solves
(\ref{plate_sol_stat}) if and only if $u$ is  a variational solution
to problem
\begin{align}\label{plate-1}
 \De^2 u +\cF(u)-p_*=C~~\mbox{in}~~\Om ,~~ u=\frac{\pd u}{\pd n}=0~~\mbox{on}~~
\pd\Omega,
\end{align}
for some constant $C$ which may depend  on $u$.
Since every variational solution to (\ref{plate-1})
is an extreme point of the functional
\[
\Psi(u)= \frac12\| \De u\|_\Om^2+ \int_\Om\Pi(u(x))dx -(p_*+C,u)_\Om,
\]
using relation (\ref{8.1.1c1}) in Assumption~\ref{A:force}
we can prove the existence of these solutions.
Thus we obtain a family of solutions to (\ref{weak_sol_stat})
parameterized by  the real parameter $C$.
To fix somehow the constant $C$ in (\ref{plate-1})
it is convenient to fix  the average of $u$.
In the case of the zero average we obtain the following
assertion.
\begin{proposition}\label{pr:stat-sol}
In addition to Assumption~\ref{A:force} we
assume that
$G_0\in V'$  and
 there exist $\eta<1/2$  and  $c\ge 0$ such that
\begin{equation}\label{u-fu}
\eta \|\De u\|^2_{\Om}+ (u, \cF(u))_\Om \ge -c,~~~\forall\,
u\in H^2_0(\Om).
\end{equation}
Then  the set $\cN_0$ of solutions $u$
to  problem (\ref{plate_sol_stat}) with the property
$\int_\Om u dx=0$ is nonempty compact set in  $\wH^2_0(\Om)$.
\end{proposition}
 \begin{proof}
Restricting  the functional $\Psi$ on   $\wH^2_0(\Om)$
we can prove the existence of its minimum point on $\wH^2_0(\Om)$.
This means that $\cN_0$ is not empty. If $u\in \wH^2_0(\Om)$
is a solution, then taking $\beta=u$ in (\ref{plate_sol_stat})
and using (\ref{u-fu}) we conclude that $\cN_0$ is bounded in
$H^2_0(\Om)$. If $\{u_n\}$ is a sequence from $\cN_0$, then
 from (\ref{plate_sol_stat}) we conclude that
\[
\|\De (u_n-u_m)\|_\Om^2\le C
\|\cF(u_n)-\cF(u_m)\|_{-1/2,\Om}\|u_n-u_m\|_{1/2,\Om}.
\]
Thus  (\ref{f-lip}) yields
$\|\De (u_n-u_m)\|_\Om\le C \|u_n-u_m\|_{2-\eps,\Om}$.
This implies that the sequence $\{u_n\}$ is relatively compact.
\end{proof}
\begin{remark}\label{re:alpha}
{\rm
A similar result can be obtain
for  the set $\cN_\al$ of solutions $u$
to  problem (\ref{plate_sol_stat}) with the property
$\langle u\rangle\equiv \int_\Om u dx'=\alpha$ with a fixed  $\al\in \R$, if
instead of \eqref{u-fu} we assume
that there exist $\eta<1/2$,  $c_\al\ge 0$ and a smooth function
$\phi$ with the property $\langle \phi\rangle=\alpha$
 such that
\begin{equation}\label{u-fu-alpha}
\eta \|\De u\|^2_{\Om}+ (u, \cF(u))_\Om-
(\phi,\cF(u))_\Om \ge -c_\al,~~~\forall\,
u\in H^2_0(\Om).
\end{equation}
Indeed, if we consider the functional $\Psi$ on
$\wH^2_{0,\al}= \left\{u\in H^2_0(\Om):\, \langle u \rangle=\alpha
\right\}$
for some fixed  constant $C$, then we
can prove the existence of a solution $u$ to (\ref{plate_sol_stat}) in
$\wH^2_{0,\al}$. Now substituting $\beta=u-\phi$
in (\ref{plate_sol_stat}) and using (\ref{u-fu-alpha})
we obtain the boundedness of the set $\cN_\al$ in $\wH^2_{0,\al}$.
To prove the compactness of $\cN_\al$ we use the same argument as in
Proposition~\ref{pr:stat-sol}.
}
\end{remark}
It follows from Proposition~\ref{pr:stat-sol} that the set
of all stationary points of $S_t$ in the space $\hch$ is nonempty
compact set and has the form
\begin{equation}\label{N-set}
\cN=\left\{ (v_*;u;0): (v_*;u)\in V_0\times \wH^2_0(\Om)~~
\mbox{solve   (\ref{v-stat}) and (\ref{plate_sol_stat})}\right\}
\end{equation}

\subsection{Asymptotical behavior} \label{sec:ab}
In this section we are interested in
global asymptotic behavior of the dynamical system  $(S_t, \hch)$.
Our main result states the existence of a compact global attractor
of finite fractal dimension.
 \par
We recall
(see, e.g., \cite{BabinVishik, Chueshov,Temam})
that
the \textit{global attractor}  of the  dynamical system  $(S_t, \hch)$
is defined as a bounded closed  set $\Ac\subset \hch$
which is  invariant ($S(t)\Ac=\Ac$ for all $t>0$) and  uniformly  attracts
all other bounded  sets:
$$
\lim_{t\to\infty} \sup\{{\rm dist}_\cH(S(t)y,\Ac):\ y\in B\} = 0
\quad\mbox{for any bounded  set $B$ in $\hch$.}
$$
The {\em fractal dimension} $\dim^X_f M$ of a compact set $M$ in a complete
metric space $X$ is defined as
\[
\dim^X_fM=\limsup_{\eps\to 0}\frac{\ln N(M,\eps)}{\ln (1/\eps)}\;,
\]
where $N(M,\eps)$ is the minimal number of closed sets in $X$ of
diameter $2\eps$ which cover~$M$.
\par
We also  recall (see, e.g., \cite{BabinVishik})
 that the \textit{unstable set} $\mathbb{M}_+(\cN)$ emanating from
some set $\cN\subset\hch$ is a subset
 of $\hch$ such that for each $z\in\mathbb{M}_+(\cN)$
there exists a full trajectory $\{y(t): t\in\R\}$ satisfying
$y(0) = z$ and ${\rm dist}(y(t),\cN) \to  0$ as $t\to -\infty$.
\par

\begin{theorem}\label{th:attractor}
Let  Assumption~\ref{A:force}  be in force.
Assume that
 $G_f(t)\equiv G_0\in V'$ is independent of $t$,
$G_{pl}(t)\equiv 0$ and (\ref{u-fu}) holds.
Then the   dynamical system  $(S_t, \widehat{\cH})$ possesses
a compact global attractor $\Ac$.
 Moreover,
\begin{enumerate}
    \item[{\bf (1)}] $\Ac=\mathbb{M}_+(\cN)$, where
 $\cN$ is the set of equilibria given by (\ref{N-set}).
\item[{\bf (2)}] This attractor has a finite  fractal dimension in $\hch$.
\item[{\bf (3)}] Any trajectory
$\gamma=\{ (v(t);u(t); u_t(t)): t\in \R\}$
from the attractor $\Ac$
possesses the properties
\begin{equation}\label{u-smth}
(v_t;u_t;u_{tt})\in L_\infty (\R;  X\times \wH^2_0(\Om)\times \Lto)
\end{equation}
and there is $R>0$ such that
\begin{equation}\label{u-smth2}
\sup_{\ga \subset \Ac}\sup_{t\in\R}\left( \|v_t\|^2_{\cO}
+\|u_t\|^2_{2,\Om}+ \|u_{tt}\|^2_{\Om}\right)\le R^2.
\end{equation}
\end{enumerate}
\end{theorem}

\begin{remark}\label{re:dyn-ch}
{\rm
We cannot state a similar result on the
 existence of a global attractor for the system $(S_t, \cH)$.
The point is that the average of $u(t)$ is preserved and
thus the system  $(S_t, \cH)$ is non-dissipative.
However  using the same procedure as for the linear case
(see Remark~\ref{re:hat-space-2})  we can  study
the long-time behavior of  $(S_t, \cH)$ by means of a family of
dissipative problems in $\hch$.
Indeed, we can decompose the solution to \eqref{fl.1}--\eqref{IC}
with the initial data $(v_0; u_0; u_1)$ into the sum $(v(t);u(t);u_t(t))=
(\bar{v}(t); \bar{u}(t); \bar{u}_t(t)) + (0;\psi;0)$,
where $\psi=(I-\widehat{P})u_0$ and $(\bar{v}(t); \bar{u}(t); \bar{u}_t(t))$
solves   \eqref{fl.1}--\eqref{fl.4}, \eqref{plBC}, \eqref{IC}
with the plate  equation
\begin{equation*}
\bar{u}_{tt} +\De^2\bar{u} +
\cF(\bar{u}+\psi)+\De^2\psi=G_{pl}(t)+p\big|_\Om
\end{equation*}
(instead of (\ref{pl_eq}))
and with the initial conditions $(v_0,\widehat{P}u_0,u_1)$.
}
\end{remark}
To obtain the result stated in Theorem~\ref{th:attractor}
it is sufficient to show
that the system is quasi-stable (in the sense of \cite{cl-book}).
For this   we use the stability properties of linear
problem (\ref{fl.1-lin})--(\ref{IC-lin}) established in
Theorem~\ref{lin_WP} to
prove the following assertion.
\begin{lemma}[Quasi-stability] \label{pr:qst}
Let $U^i(t)=(v^i(t);u^i(t);u^i_t(t))$, $i=1,2$,
 be two weak solutions with
initial data $U^i_0=(v^i_0;u^i_0;u^i_1)$ from  $\hch$
such that
 $\|U_0^i\|_\cH\le R$, $i=1,2$, then
their difference
\[
Z(t)=U^1(t)-U^2(t)\equiv (v(t);u(t);u_t(t))
\]
satisfies the relation
\begin{equation}\label{q-stab}
\|Z(t)\|^2_\cH\le M_R e^{-\ga_* t} \|Z_0\|^2_\cH+M_R \int_0^t
e^{-\ga_*(t-\tau)}
\|u(\tau)\|^2_{\Om} d\tau
\end{equation}
for some positive constant $M_R$ and $\ga_*$.
\end{lemma}
\begin{proof}
We consider $(v(t);u(t))$ as a solution to
to linear problem (\ref{fl.1-lin})--(\ref{plBC-lin}) with
$G_f\equiv 0$ and $G_{pl}(t)=-\cF(u^1(t))+ \cF(u^2(t))$.
Therefore
it follows from (\ref{f-lip})  and  (\ref{exp-st-nhm})
that
\[
\|Z(t)\|^2_\cH\le M e^{-\ga t} \|Z_0\|^2_\cH+ C_R \int_0^t
e^{-\ga(t-\tau)}
\|u(\tau)\|^2_{2-\eps, \Om} d\tau.
\]
Hence the
interpolation relation
\[
\|u\|^2_{2-\eps,\Om}\le \eta\|Z\|^2_\cH +  c_\eta \|u\|^2_{\Om},
~~~\forall\, \eta>0,
\]
via Gronwall's type argument, implies the conclusion in (\ref{q-stab}).
\end{proof}
\subsubsection*{Proof  of Theorem~\ref{th:attractor}}
Lemma~\ref{pr:qst} means that the dynamical system
 $(S_t, \widehat{\cH})$ is quasi-stable in the sense of Definition 7.9.2~\cite{cl-book}.
Therefore by Proposition~7.9.4~\cite{cl-book}
 $(S_t, \widehat{\cH})$ is asymptotically smooth.
Since the system is gradient,
the boundedness of
 the set of the stationary points
implies that there exists a compact global
attractor.
Moreover,
the standard results on gradient systems with compact attractors
 (see, e.g., \cite{BabinVishik,Chueshov,Temam}) give us that $\Ac=\mathbb{M}_+(\cN)$.
\par
Since  $(S_t, \widehat{\cH})$
is quasi-stable   the finiteness of fractal dimension ${\rm dim}_f\Ac$
follows from  Theorem~7.9.6~\cite{cl-book}.
\par
 To obtain the result on regularity stated in (\ref{u-smth}) and
(\ref{u-smth2}) we apply  Theorem~7.9.8~\cite{cl-book}.

\appendix

\section{Appendix: Generator of  linear semigroup}\label{ap:generator}
 To find
the structure of  the semigroup $T_t$ generated by
(\ref{fl.1-lin})--(\ref{IC-lin}) in the space $\hch$
 we note that
the evolution problem in (\ref{weak_sol_d2}) with  $G_f\equiv 0$ and $G_{pl}\equiv0$ can be written in the form
\[
\frac{d}{dt}\left[(v,\psi)_\cO+ (u(t),\chi)_\Om+(w(t),\beta)_\Om  \right] +\cA(U(t),\Psi)=0,
\]
where $U=(v;u;w)$ is an element from $C(\R_+;\hch)$ with $v\in L_2^{loc}(\R_+; V)$ and
$v\big|_\Om =(0;0;w)$. The text function $\Psi=(\psi;\chi;\beta)$ belongs to the space
\[
 \cW   \equiv \left\{(\psi;\chi;\beta)\in 
W\times \wH^2_0(\Om)\times \Hto\, :\; \psi\big|_\Om =(0;0;\beta)\right\}\subset\hch,
\]
and the  bilinear form  $\cA(U,\Psi)$ is
defined  by the relation
\[
\cA(U,\Psi)= \nu(\g v,\g\psi)_{\cO} -(w, \chi)_\Om+ (\De u, \De \beta)_\Om.
\]
Thus to describe  the domain of the  generator we need
 to describe all elements $U=(v;u;w)$ from 
 \[
\cV \equiv\left\{(v;u;w)\in V\times \wH^2_0(\Om)\times \Lto\, :\; v\big|_\Om =(0;0;w)\right\}\subset\hch 
\]
which solve
the variational equation of the form
\[
\cA(U,\Psi)=(v,\psi)_\cO+ (u,\chi)_\Om+(w,\beta)_\Om ,~~~ \forall\,\Psi=(\psi;\chi;\beta)\in\cW,
\]
where $F=(f_0;f_1;f_2)$ is a given element from $\hch$.
Taking $\psi\equiv 0$ one can see that $f_1=-w\in \wH^2_0(\Om)$.
Therefore we arrive at the relation
\begin{equation}\label{var-g}
  \nu(\g v,\g\psi)_{\cO} + (\De u, \De \beta)_\Omega=
 (f_0,\psi)_{\cO} + (f_2, \beta)_\Om.
\end{equation}
By Proposition~\ref{pr:stokes} we have that $N_0w\in V\cap [H^2(\cO)]^3$
and the corresponding pressure $p_w$ (defined by (\ref{fl.n0}))
belongs to the class $H^1(\cO)/\R$.
Since
\[
\nu(\g N_0w,\g \psi)_\cO=-\nu(\Delta N_0w,\psi)_\cO=-(\g p_w,\psi)_\cO=
-(p_w,\beta)_\Om.
\]
We can rewrite (\ref{var-g}) in the form
\begin{equation}\label{var-g1}
  \nu(\g [v-N_0w],\g\psi)_{\cO} + (\De u, \De \beta)_\Omega=
 (f_0,\psi)_{\cO} + (f_2+p_w, \beta)_\Om
\end{equation}
for any $\psi\in W$.
If we take now $\psi\in \widetilde{V}=\{ v\in V :\, v|_{\pd\cO}=0\}$, then
we obtain that $\tilde{v}=v-N_0w\in \widetilde{V}$ solve the problem
\[
 -\nu\Delta \tilde{v}+\nabla p= f_0, ~~
   \di \tilde{v}=0 ~~ {\rm in} ~ \cO;
\quad \tilde{v}=0 ~~ {\rm on}~\pd\cO.
\]
Since $f_0\in X$, this implies that
$\tilde{v}\in \widetilde{V}\cap [H^2(\cO)]^3$
and thus $v\in V\cap [H^2(\cO)]^3$.
Therefore from (\ref{var-g}) we have that
\begin{equation}\label{var-g2}
  (P_S[-\nu \Delta v]- f_0,\psi)_{\cO} +
(\De u,   \De \beta)_\Omega = (f_2,  \beta)_\Omega
\end{equation}
for every $\psi\in X$ with $\psi^3|_\Om=\beta\in\Hto$,
where $P_S$ is the orthoprojector in   $[L_2(\cO)]^3$
on $X$. This implies that
\[
P_S[-\nu \Delta v]- f_0 ~\bot~ \widetilde{X}=\{ u\in X :\, (u, n)=0~\mbox{on}~\pd\cO\}
\]
Therefore  (see, e.g., (2.70) in \cite{AHKM03})
there exists $q\in H^1(\cO)$ such that
\begin{equation}\label{var-g3}
  P_S[-\nu \Delta v]- f_0 =- \g q,~~
\De q=0~\mbox{in}~\cO,~~ \frac{\pd q}{\pd n}\Big|_S=0.
\end{equation}
Substitution in (\ref{var-g2}) yields
$(\De u, \De \beta)_\Omega= (f_2+ q,  \beta)_\Omega$
which implies that $u\in (H^4\cap H^2_0)(\Om)$.
On the other hand, if we take $\psi=N_0\beta$ in (\ref{var-g1}), then
due to the relation $(\g [v-N_0w],\g N_0\beta)_{\cO}=0$
 we obtain
\[
(\De u, \De \beta)_\Omega= (N^*_0f_0 +f_2+ p_w,  \beta)_\Omega,
~~\forall\, \beta\in \Hto.
\]
Thus,  since the function $q$ is defined up to a constant,
we can suppose  that
\begin{equation}\label{var-g5}
 q\big|_\Om = N^*_0f_0+ p_w= -f_2+
\De^2 u -\int_\Om \De^2 ud x' \in \Lto.
\end{equation}
Let us denote by $\cG :H^{1/2}_*(\Om)\mapsto X$ the mapping $r\mapsto \g q$,
where $q\in H^1(\cO)$ solve the problem
\[
\De q=0~\mbox{in}~\cO,~~ \frac{\pd q}{\pd n}\Big|_S=0,~~ q\big|_\Om =r.
\]
Let $$\bar{X}=\{u\in X :\, \ga_nu\equiv (u,n)\big|_\Om\in L_2(\Om)\}$$
equipped with the graph norm
$\|u\|_{\bar{X}}^2=  \|u\|_{\cO}^2+\|\ga_n u\|_{\Om}^2$.
It is obvious that the trace operator $\ga_n$ is bounded from
$\bar{X}$ into $L^2(\Om)$.
One can see from calculations on smooth functions that
\[
(\cG r,\psi)_\cO =\int_{\pd\cO} q (\psi,n)dS=
\int_{\Om} r (\psi,n)dx',
~~\forall\, \psi\in \bar{X}.
\]
and therefore
\begin{equation}\label{g-positive}
(\cG \ga_n\phi ,\psi)_\cO =( \ga_n\phi, \ga_n\psi)_\Om,
~~\forall\, \phi, \psi\in \bar{X}.
\end{equation}
 Consequently  the operator $\Gamma=\cG\ga_n$ 
 can be extended to a bounded operator on $\bar{X}$.
Moreover, by \eqref{g-positive} $\Gamma$ is nonnegative. With this  
 operator $\Gamma$ using the fact that
that $f_2=\ga_n f_0$ we can write
(\ref{var-g3})  in the form
\[
f_0+\Gamma f_0=  P_S[-\nu \Delta v] +\cG\left[
\De^2 u -\int_\Om \De^2 ud x'\right]\equiv  \sQ(v,u).
\]
This leads to the following description
of the generator  $\cA$:
\[
\sD(\cA)=\left\{ (v;u;w)\in\hch\, \left|\begin{array}{l}
v\in V\cap [H^2(\cO)]^3,   u\in H^4(\Om), \\ w\in \Hto,
\ga_n[\sQ(v,u)]\in L_2(\Om)
\end{array}\right.
\right\}
\]
and
\[
\cA\begin{bmatrix} v\\ u\\ w\end{bmatrix} =
 \begin{bmatrix} (1+\Gamma)^{-1}\sQ(v,u) \\ -w \\
\De^2 u -\int_\Om \De^2 ud x'-p_w
-N_0^*(1+\Gamma)^{-1}\sQ(v,u)
\end{bmatrix}.
\]
We can also write the operator $\cA$ in the form
\[
\cA\begin{bmatrix} v\\ u\\ w\end{bmatrix} =
 \begin{bmatrix} (1+\Gamma)^{-1}\sQ(v,u) \\ -w \\
\ga_n (1+\Gamma)^{-1}\sQ(v,u)
\end{bmatrix}.
\]

\end{document}